\documentclass[11pt]{amsart}
\usepackage{amsopn}
\usepackage{amssymb, amscd}
\usepackage{multirow}
\usepackage{graphicx, graphics, epsfig}
\usepackage{faktor}
\usepackage{enumerate}
\usepackage{tikz}
\usepackage{pgfplots}
\usepackage{hyperref}
\usepackage{color}

\newcommand{\nc}{\newcommand}

\nc{\fg}{\mathfrak{f} } \nc{\vg}{\mathfrak{v} } \nc{\wg}{\mathfrak{w} }
\nc{\zg}{\mathfrak{z} } \nc{\ngo}{\mathfrak{n} } \nc{\kg}{\mathfrak{k} }
\nc{\mg}{\mathfrak{m} } \nc{\bg}{\mathfrak{b} } \nc{\ggo}{\mathfrak{g} } \nc{\eg}{\mathfrak{e} }
\nc{\ggob}{\overline{\mathfrak{g}} } \nc{\sog}{\mathfrak{so} }
\nc{\sug}{\mathfrak{su} } \nc{\spg}{\mathfrak{sp} } \nc{\slg}{\mathfrak{sl} }
\nc{\glg}{\mathfrak{gl} } \nc{\cg}{\mathfrak{c} } \nc{\rg}{\mathfrak{r} }
\nc{\hg}{\mathfrak{h} } \nc{\tg}{\mathfrak{t} } \nc{\ug}{\mathfrak{u} }
\nc{\dg}{\mathfrak{d} } \nc{\ag}{\mathfrak{a} } \nc{\pg}{\mathfrak{p} }
\nc{\sg}{\mathfrak{s} } \nc{\affg}{\mathfrak{aff} } \nc{\qg}{\mathfrak{q} } \nc{\lgo}{\mathfrak{l} }

\nc{\pca}{\mathcal{P}} \nc{\nca}{\mathcal{N}} \nc{\lca}{\mathcal{L}}
\nc{\oca}{\mathcal{O}} \nc{\mca}{\mathcal{M}} \nc{\tca}{\mathcal{T}}
\nc{\aca}{\mathcal{A}} \nc{\cca}{\mathcal{C}} \nc{\gca}{\mathcal{G}}
\nc{\sca}{\mathcal{S}} \nc{\hca}{\mathcal{H}} \nc{\bca}{\mathcal{B}}
\nc{\dca}{\mathcal{D}} \nc{\eca}{\mathcal{E}} \nc{\wca}{\mathcal{W}}

\nc{\vp}{\varphi} \nc{\ddt}{\tfrac{d}{dt}} \nc{\dsdt}{\tfrac{d^2}{dt^2}} \nc{\dds}{\tfrac{d}{ds}}
\nc{\dpar}{\tfrac{\partial}{\partial t}} \nc{\im}{\mathrm{i}}

\nc{\SO}{\mathrm{SO}} \nc{\Spe}{\mathrm{Sp}} \nc{\Sl}{\mathrm{SL}}
\nc{\SU}{\mathrm{SU}} \nc{\Or}{\mathrm{O}} \nc{\U}{\mathrm{U}} \nc{\Gl}{\mathrm{GL}}
\nc{\Se}{\mathrm{S}} \nc{\Cl}{\mathrm{Cl}} \nc{\Spin}{\mathrm{Spin}}
\nc{\Pin}{\mathrm{Pin}} \nc{\G}{\mathrm{GL}_n(\RR)} \nc{\g}{\mathfrak{gl}_n(\RR)}

\nc{\RR}{{\Bbb R}} \nc{\HH}{{\Bbb H}} \nc{\CC}{{\Bbb C}} \nc{\ZZ}{{\Bbb Z}}
\nc{\FF}{{\Bbb F}} \nc{\NN}{{\Bbb N}} \nc{\QQ}{{\Bbb Q}} \nc{\PP}{{\Bbb P}} \nc{\OO}{{\Bbb O}}

\nc{\vs}{\vspace{.2cm}} \nc{\vsp}{\vspace{1cm}} \nc{\ip}{\langle\cdot,\cdot\rangle}
\nc{\ipp}{(\cdot,\cdot)} \nc{\la}{\langle} \nc{\ra}{\rangle} \nc{\unm}{\tfrac{1}{2}}
\nc{\unc}{\tfrac{1}{4}} \nc{\und}{\tfrac{1}{16}} \nc{\no}{\vs\noindent}
\nc{\lam}{\Lambda^2(\RR^n)^*\otimes\RR^n} \nc{\tangz}{{\rm T}^{\rm Zar}}
\nc{\nor}{{\sf n}}  \nc{\mum}{/\!\!/} \nc{\kir}{/\!\!/\!\!/}
\nc{\Ri}{\tfrac{4\Ric_{\mu}}{||\mu||^2}} \nc{\ds}{\displaystyle}
\nc{\ben}{\begin{enumerate}} \nc{\een}{\end{enumerate}} \nc{\f}{\tfrac}
\nc{\lb}{[\cdot,\cdot]} \nc{\isn}{\tfrac{1}{||v||^2}}
\nc{\gkp}{(\ggo=\kg\oplus\pg,\ip)} \nc{\ukh}{(\ug=\kg\oplus\hg,\ip)}
\nc{\tgkp}{(\tilde{\ggo}=\kg\oplus\pg,\ip)}
\nc{\wt}{\widetilde}
\nc{\iop}{\mathtt{i}} \nc{\jop}{\mathtt{j}}

\nc{\Hess}{\operatorname{Hess}} \nc{\ad}{\operatorname{ad}}
\nc{\Ad}{\operatorname{Ad}} \nc{\rank}{\operatorname{rk}}
\nc{\Irr}{\operatorname{Irr}} \nc{\End}{\operatorname{End}}
\nc{\Aut}{\operatorname{Aut}} \nc{\Inn}{\operatorname{Inn}}
\nc{\Der}{\operatorname{Der}} \nc{\Ker}{\operatorname{Ker}}
\nc{\Iso}{\operatorname{Iso}} \nc{\Diff}{\operatorname{Diff}}
\nc{\Lie}{\operatorname{L}} \nc{\tr}{\operatorname{tr}} \nc{\dif}{\operatorname{d}}
\nc{\sen}{\operatorname{sen}} \nc{\modu}{\operatorname{mod}}
\nc{\CRic}{\operatorname{PP}} \nc{\Cric}{\operatorname{P}} \nc{\Ricci}{\operatorname{Ric}}
\nc{\sym}{\operatorname{sym}} \nc{\herm}{\operatorname{herm}} \nc{\symac}{\operatorname{sym^{ac}}}
\nc{\symc}{\operatorname{sym^{c}}} \nc{\scalar}{\operatorname{Sc}}
\nc{\grad}{\operatorname{grad}} \nc{\ricci}{\operatorname{Rc}} \nc{\kil}{\operatorname{B}} \nc{\cas}{\operatorname{C}} \nc{\lic}{\operatorname{L}}
\nc{\Nor}{\operatorname{Norm}}  \nc{\ricc}{\operatorname{Rc^{c}}}
\nc{\Ricc}{\operatorname{Ric^{c}}} \nc{\ricac}{\operatorname{Rc^{ac}}}
\nc{\Ricac}{\operatorname{Ric^{ac}}} \nc{\Riem}{\operatorname{Rm}} \nc{\Sec}{\operatorname{Sec}}
\nc{\riccig}{\operatorname{ric^{\gamma}}} \nc{\mm}{\operatorname{m}} \nc{\Mm}{\operatorname{M}}
\nc{\Le}{\operatorname{L}} \nc{\tang}{\operatorname{T}}
\nc{\level}{\operatorname{level}} \nc{\rad}{\operatorname{r}}
\nc{\abel}{\operatorname{ab}} \nc{\CH}{\operatorname{CH}} \nc{\Cone}{{\mathcal C}} \nc{\CCone}{\operatorname{CC}} \nc{\CP}{{\mathcal P}}
\nc{\mcc}{\operatorname{mcc}} \nc{\Adj}{\operatorname{Adj}}
\nc{\Order}{\operatorname{O}}  \nc{\inj}{\operatorname{inj}} \nc{\proy}{\operatorname{pr}}
\nc{\vol}{\operatorname{vol}} \nc{\Diag}{\operatorname{Dg}} \nc{\Diagg}{\operatorname{Diag}}
\nc{\Spec}{\operatorname{Spec}} \nc{\Ima}{\operatorname{Im}} \nc{\Rea}{\operatorname{Re}}
\nc{\spann}{\operatorname{span}} \nc{\Aff}{\operatorname{Aff}} \nc{\E}{\operatorname{E}} \nc{\id}{\operatorname{id}} \nc{\dete}{\operatorname{det}} \nc{\Crit}{\operatorname{Crit}} \nc{\val}{\operatorname{val}}

\theoremstyle{plain}
\newtheorem{theorem}{Theorem}[section]

\theoremstyle{definition}
\newtheorem{definition}[theorem]{Definition}

\theoremstyle{remark}
\newtheorem{remark}[theorem]{Remark}

\title{Homogeneous Einstein metrics and local maxima of the Hilbert action}

\author{Jorge Lauret}  \author{Cynthia E. Will}

\address{FaMAF, Universidad Nacional de C\'ordoba and CIEM, CONICET (Argentina)}
\email{jorgelauret@unc.edu.ar}  \email{cynthia.will@unc.edu.ar}

\thanks{This research was partially supported by grants from FONCYT and Univ. Nac. de C\'ordoba   (Argentina).}

\date{\today}

\begin{document}

\maketitle

\begin{abstract}
In this short note, three infinite families of neutrally stable homogeneous Einstein metrics are ruled out as candidates for local maxima of the Hilbert action.  
\end{abstract}

\tableofcontents

%
%\iffalse 
%

%
%\fi
%

\section{Introduction}\label{intro}

As well known, Einstein metrics on a compact differentiable manifold $M$ are precisely the critical points of the total scalar curvature functional (or Hilbert action)
\begin{equation}\label{sct}
\widetilde{\scalar}(g):=\int_M \scalar(g)\; d\vol_g,
\end{equation}
restricted to the space $\mca_1$ of all Riemannian metrics on $M$ of some fixed volume (see \cite[4.21]{Bss}).  Any Einstein metric is a saddle point, but mainly due to the fact that the Hessian is positive on conformal variations;  indeed, when $\widetilde{\scalar}$ is further restricted to the submanifold
$$
\cca_1:=\{ g\in\mca_1:\scalar(g)\,\mbox{is a constant function on}\, M\}, 
$$ 
the nullity and coindex of critical points are both finite and hence the existence of local maxima of $\widetilde{\scalar}|_{\cca_1}$ emerges as a possibility.  More precisely, if $g\in\mca_1$ is Einstein, then
\begin{equation}\label{Cdec}
T_g\cca_1 = \lca_{\mathfrak{X}(M)}g \oplus^{\perp_g} \tca\tca_g, 
\end{equation}
where $\lca_{\mathfrak{X}(M)}g=\Ima \delta_g^*=T_g\Diff(M)\cdot g$ is the space of trivial variations and $\tca\tca_g:=\Ker\delta_g\cap\Ker\tr_g$ is the subspace of divergence-free (or transversal) and traceless symmetric $2$-tensors, so-called {\it TT-tensors} (see \cite[4.44-4.46]{Bss}).  $\lca_{\mathfrak{X}(M)}g$ is therefore contained in the kernel of the Hessian $\widetilde{\scalar}''_g$ of $\widetilde{\scalar}$ at $g$ and it was independently proved by Koiso and Berger that $\widetilde{\scalar}''_g|_{\tca\tca_g}$ is negative definite on the orthogonal complement of a (possibly trivial) finite-dimensional vector subspace of $\tca\tca_g$ (see \cite[4.60]{Bss}).  

\begin{figure}
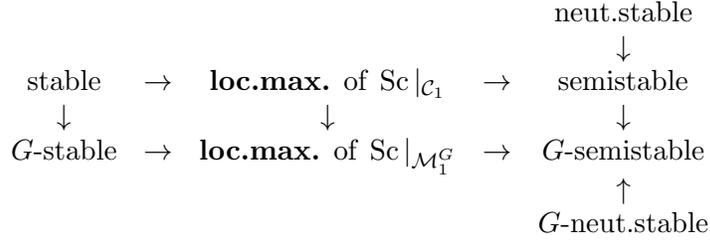

$$
\begin{array}{ccccc}
 && &  &  \text{neut.stable}\\ 
& &&& \downarrow\\ 
%
%\text{RF dyn.stable}& \leftarrow & 
\text{stable} & \rightarrow & \text{{\bf loc.max.} of}\; \scalar|_{\cca_1} & \rightarrow & \text{semistable} \\ 
 \downarrow && \downarrow && \downarrow \\ 
 G\text{-stable} & \rightarrow & \text{{\bf loc.max.} of}\; \scalar|_{\mca_1^G} & \rightarrow &  G\text{-semistable}\\ 
& &&& \uparrow \\ 
&& &  & G\text{-neut.stable} \\ 
\end{array}
$$
\caption{Some stability types for an Einstein metric.}\label{stabtypes}
\end{figure}

\begin{definition}\label{stab-def-2} (see Figure \ref{stabtypes}).  
An Einstein metric $g\in\mca_1$ is said to be,  
\begin{enumerate}[{\small $\bullet$}] 
\item {\it stable} (or linearly stable): $\widetilde{\scalar}''_g|_{\tca\tca_g}<0$.  In particular, $g$ is a local maximum of $\widetilde{\scalar}|_{\cca_1}$ since by \eqref{Cdec}, $\tca\tca_g$ exponentiates into a slice for the $\Diff(M)$-action (see \cite[12.22]{Bss} or \cite[Lemma 2.6.3]{Krn2}).  

\item {\it semistable}: $\widetilde{\scalar}''_g|_{\tca\tca_g}\leq 0$ ; and otherwise {\it unstable}.  Note that semistability must hold for any local maximum of $\widetilde{\scalar}|_{\mca_1}$. 

\item {\it neutrally stable}: $\widetilde{\scalar}''_g|_{\tca\tca_g}\leq 0$ and has nonzero kernel, i.e., it is in addition {\it infinitesimally deformable}.  In this case, $g$ may or may not be a local maximum of $\widetilde{\scalar}|_{\mca_1}$.   Note that any semistable Einstein metric is either stable or neutrally stable.  
 
\item {\it dynamically stable}: for any metric $g_0$ near $g$, the normalized Ricci flow starting at $g_0$ exists for all $t\geq 0$ and converges modulo diffeomorphisms, as $t\to\infty$, to an Einstein metric near $g$.  Stability is a necessary condition for $\nu$-stability, where $\nu$ is the Perelman's entropy functional (see \cite{CaoHe,WngWng}), which is in turn a sufficient condition for Ricci flow dynamical stability (see \cite[Theorem 1.3]{Krn} or \cite[Corollary 6.2.5]{Krn2}).       
\end{enumerate}
\end{definition}

Among compact irreducible symmetric spaces, it is well known that most of them are stable, some of them are unstable and precisely the following ones are neutrally stable (see \cite{Kso, GsqGld, SmmWng, Sch}):   
\begin{equation}\label{nsss}
\begin{array}{c}
\SU(n), \quad \SU(n)/ \SO(n), \quad \SU(2n)/ \Spe(n),\quad n \geq 3, \\ 
\SU(p + q)/\Se(\U(p) \times \U(q)),\quad p \geq q \geq 2 \quad\mbox{and} \quad E_6/F_4.  
\end{array}
\end{equation}
Curiously enough, stable irreducible symmetric spaces are the only known local maxima of $\widetilde{\scalar}|_{\cca_1}$ with $\scalar>0$.  Moreover, the only neutrally stable cases which have been proved not to be local maxima are $\SU(3)$ (see \cite{Jns}) and $\SU(2n)/ \Spe(n)$ (see \cite[Example 6.7]{BhmWngZll}).  In both cases $G$-invariant variations for some Lie group $G$ acting transitively on $M$ were used.  

In fact, in the $G$-invariant setting, it is natural to study the $G$-invariant versions of the notions defined above, i.e., $G$-stability, etc., by restricting the Hilbert action to the finite-dimensional submanifold of $G$-invariant metrics $\mca_1^G$ (see Figure \ref{stabtypes}).  We are therefore lead to study the Hessian $\widetilde{\scalar}''_g$ restricted to the finite-dimensional subpace of $G$-invariant $TT$-tensors $\tca\tca_g^G$ (see \cite{stab-tres, stab-dos, stab}).  Note that $\mca_1^G$ is a singleton for any irreducible symmetric space in \eqref{nsss} relative to the given transitive groups, as they are all isotropy irreducible homogeneous spaces, so one needs to consider a smaller transitive group in order to enlarge $\mca_1^G$ and get $G$-invariant variations.  From the list \eqref{nsss}, this can only be done for $\SU(2n)/\Spe(n)$ and $\SU(n)$.  We use these new presentations to find an explicit curve $g(t)\in\mca_1^G$ through the neutrally stable symmetric metric $g_0$ such that $\scalar(g(t))$ has an inflection point, which implies that $g_0$ is not a local maximum of $\scalar|_{\mca_1^G}$ and so neither of $\widetilde{\scalar}|_{\cca_1}$.  The formulas for the Hessian of $\scalar$ given in \cite{stab-tres, stab-dos} were very helpful to find the curve $g(t)$.  

\begin{theorem}\label{main}
The following Einstein metrics are not local maxima of $\widetilde{\scalar}|_{\cca_1}$: 
\begin{enumerate}[{\rm (i)}]
\item The bi-invariant metric on $\SU(n)$ for any $n\geq 3$, which is not a local maximum of $\scalar$ among left-invariant metrics of volume one.  

\item \cite[Example 6.8]{BhmWngZll} The symmetric metric on $\SU(2n)/\Spe(n)$ for any $n\geq 3$,  which is not a local maximum of $\scalar$ among $\SU(2n-1)$-invariant metrics of volume one.    
\end{enumerate}
\end{theorem}

\begin{remark}
Part (i) is stated in \cite[Corollary 5.4]{BhmWngZll}, but the proof refers to \cite{Jns}, where it is done only for $n=3$.  It was proved in \cite{CaoHe} that this metric is not a local maximum of the Perelman $\nu$-entropy and so it is dynamically unstable by using conformal variations (see also \cite[Theorem E]{BttHllMrpWld}).  Concerning part (ii), we give here an elementary alternative proof.   
\end{remark}

\begin{remark}
After the first version of this paper appeared in arXiv, Yu.G.\ Nikonorov uploaded the English translation \cite{Nkn} of a paper in Russian which is part of his dissertation \cite{Nkn2}, where part (i) is proved (see \cite[Theorem 3, 4)]{Nkn}).  
\end{remark}

On the other hand, only few examples of $G$-stable or local maxima of $\scalar|_{\mca_1^G}$ have been found in the explorations made in \cite{stab-tres, stab-dos, stab} throughout standard metrics, generalized Wallach spaces and flag manifolds with $b_2(M)=1$.  Beyond the case of homogeneous spaces with two isotropy summands (see \S\ref{r2-sec}), they are all either standard or K\"ahler metrics.  $G$-neutrally stability is even more scarce, only the following examples were found: the standard metrics on $\SO(2n)/T^n$, $n\geq 4$ (see \cite{stab}).  Our second result rules out these $G$-neutrally stable Einstein metrics as candidates for local maxima of $\widetilde{\scalar}|_{\cca_1}$.    

\begin{theorem}\label{main2}
The standard metric on the full flag manifold $\SO(2n)/T^n$ is not a local maximum of $\scalar|_{\mca_1^G}$ (nor of $\widetilde{\scalar}|_{\cca_1}$) for any $n\geq 4$.  
\end{theorem}

\begin{remark}
It follows from \cite[Theorem 5.1]{BhmWngZll} that none of the Einstein metrics in Theorems \ref{main} and \ref{main2} realizes the Yamabe invariant of $M$ (i.e., $\scalar(g)$ is not the supreme among all Yamabe metrics of $M$, which are those with the smallest scalar curvature in its unit volume conformal class). 
\end{remark}

In particular, the only known homogeneous Einstein metrics with $\scalar>0$ which still have a chance to be new examples of local maxima of the Hilbert action $\widetilde{\scalar}|_{\cca_1}$ are:

\begin{enumerate}[{\small $\bullet$}] 
\item the irreducible symmetric spaces listed in \eqref{nsss} different from $\SU(2n)/\Spe(n)$ and $\SU(n)$, 

\item the remaining isotropy irreducible homogeneous spaces, 

\item the local maxima of $\scalar|_{\mca_1^G}$ in the two isotropy summands case (see \S\ref{r2-sec}), 

\item and the $G$-stable Einstein metrics found in \cite{stab-dos} and \cite{stab}.    
\end{enumerate}

\begin{remark}\label{new}
After the first version of the present paper was uploaded to arXiv, the $G$-stable Einstein metric on $E_7/\SO(8)$ (see \cite{stab-dos}) has been proved to be stable in \cite{SchSmmWng}.  
\end{remark}

\vs \noindent {\it Acknowledgements.}  We are very grateful with Christopher B\"ohm and Wolfgang Ziller for many helpful conversations.

\section{Preliminaries}\label{preli} 

Let $M$ be a compact connected differentiable manifold of dimension $d$.  We assume that $M$ is homogeneous and fix an almost-effective transitive action of a compact Lie group $G$ on $M$, which determines a presentation $M=G/K$ of $M$ as a homogeneous space, where $K\subset G$ is the isotropy subgroup at some point $o\in M$.  Let $\mca^G\subset\sca^2(M)^G$ denote the open cone of all $G$-invariant Riemannian metrics on $M$, where $\sca^2(M)^G$ is the finite-dimensional vector space of all $G$-invariant symmetric $2$-tensors.  Note that $\mca^G$ is therefore a differentiable manifold of dimension between $1$ and $\tfrac{d(d+1)}{2}$.  

It is well known that $g\in\mca_1^G$ is Einstein if and only if $g$ is a critical point of the scalar curvature functional $\scalar:\mca_1^G\longrightarrow \RR$, where $\mca^G_1\subset\mca^G$ is the codimension one submanifold of all metrics with a given volume.  

We fix a bi-invariant inner product $Q$ on $\ggo$ and consider the $Q$-orthogonal reductive decomposition $\ggo=\kg\oplus\pg$, and as a background metric, $Q|_{\pg}\in\mca^G$.  For each $Q$-orthogonal decomposition,
\begin{equation}\label{dec2}
\pg=\pg_1\oplus\dots\oplus\pg_r,
\end{equation}
in $\Ad(K)$-invariant  subspaces (not necessarily $\Ad(K)$-irreducible), we denote the corresponding structural constants by
\begin{equation}\label{ijk}
[ijk]:=\sum_{\alpha,\beta,\gamma} Q([e_{\alpha}^i,e_{\beta}^j], e_{\gamma}^k)^2,
\end{equation}
where $\{ e_\alpha^i\}$, $\{ e_{\beta}^j\}$ and $\{ e_{\gamma}^k\}$ are $Q$-orthonormal basis of $\pg_i$, $\pg_j$ and $\pg_k$, respectively.  Note that $[ijk]$ is invariant under any permutation of $ijk$. 

For each $g\in\mca^G$, there exists at least one decomposition as in \eqref{dec2} which is also $g$-orthogonal and so
\begin{equation}\label{metric}
g=x_1Q|_{\pg_1}+\dots+x_rQ|_{\pg_r}, \qquad x_i>0,
\end{equation}
and the scalar curvature of $g=(x_1,\dots,x_r)$ is given by 
\begin{equation}\label{scal}
\scalar(g) = \unm\sum_k\tfrac{b_kd_k}{x_k}-\tfrac{1}{4}\sum_{i,j,k} \tfrac{x_k}{x_i x_j}[ijk].  
\end{equation}
Note that if we set $Q:=-\kil_\ggo$ in the case when $G$ is semisimple, where $\kil_\ggo$ is the Killing form of $\ggo$, then $b_k=1$ for all $k$.

\section{Two isotropy summands}\label{r2-sec}

For homogeneous spaces with only two isotropy summands, the above formulas considerably simplify.   We refer to \cite{DckKrr, He} for a complete classification of these spaces ($43$ infinite families and $78$ isolated examples) together with their $G$-invariant Einstein metrics (only $19$ isolated examples do not admit any).  We assume that the summands are inequivalent since, according to \cite{DckKrr}, the only case where this is not the case is $S^7\times S^7=\SO(8)/G_2$, which admits exactly two non-homothetic Einstein metrics, the standard and the product metric (see \cite{Krr}), both known to be $G$-unstable (see \cite[Table 1, 9]{stab} for the standard metric).  Thus $\dim{\mca_1^G}=1$ and so $\scalar|_{\mca_1^G}$ is a one-variable function. 

A summary of the possible behaviors in this case concerning the existence of local maxima of $\scalar|_{\mca_1^G}$ and $G$-neutrally stable Einstein metrics follows:  

\begin{enumerate}[{\rm (i)}] 
\item If $K$ is a maximal subgroup of $G$ ($5$ infinite families and $13$ isolated examples, see \cite[Section 6]{DckKrr}), then there exists a global maximum of $\scalar|_{\mca_1^G}$ and at most two more Einstein metrics (see \cite[Section 3.2]{Krr}).  Except for a few exceptions (see \cite[Section 6]{stab-tres}), a complete picture is not available in the literature for these spaces.  

\item Assume in what follows that $K$ is not a maximal subgroup of $G$.  There are at most two Einstein metrics and global minima or maxima are not allowed (see \cite[Theorem 1.3]{WngZll} and \cite[Section 3]{DckKrr}.  

\item In the case when the space admits a unique invariant Einstein metric, such a metric is necessarily an inflection point of $\scalar|_{\mca_1^G}$ as global minima and maxima do not exist, and so it is never a local maximum of $\scalar|_{\mca^G_1}$.  Note that it is $G$-neutrally stable.  The list of such spaces, including the intermediate subgroup $K< H< G$, can be extracted from the lists obtained in \cite{DckKrr} and it is given by: 

\begin{enumerate}[{\tiny $\bullet$}] 
  \item I.24: $\SU(k) < \U(k) < \SO(2k)$, $k \ge 3$.

  \item I.25: $\SO(3)\U(1) < \U(3) < \SO(6)$.

  \item II.2: $\Se(\U(2)\U(1)) < S(\U(3)\U(1)) < \SU(4)$.

  \item II.6: $\Se(\SO(n)\times \U(1) \times \U(m)) < \Se(\U(n)\U(m)) < \SU(n+m)$, $n=m^2+2$, $m\ge 1$.

  \item II.7:  $\SU(m)\times \SU(n) < \Se(\U(n)\U(m)) < \SU(n+m)$, $n,m\ge 2$.

  \item III.8: $\SU(k) < \U(k) < \Spe(k)$, $k \ge 3$.

  \item IV.7: $\Spin(10) < \Spin(10) \SO(2) < E_6$.

 \item IV.16: $\SO(6) \SU(2) < \SU(6)\SU(2) < E_6$ (standard).

 \item IV.19: $E_6 < \E_6 \times \SO(2) < E_7$.
\end{enumerate}

\item As far as we know, the existence of Einstein metrics on the $5$ spaces added in \cite[Appendix A]{He} has not been studied yet.  We note that the space $E_8/\Spin(9)$ was missed in the papers \cite{DckKrr, He}.  
\end{enumerate}

The only space in part (iii) for which the standard metric $g_{\kil}$ is Einstein is $M=E_6/\SU(2)\times \SO(6)$ (see \cite{WngZll2}).  In what follows, we give an alternative proof of the fact that $g_{\kil}$ is not a local maximum of $\scalar|_{\mca^G_1}$.  

There is a reductive decomposition $\ggo=\kg\oplus\pg$ and a decomposition $\pg=\pg_1\oplus\pg_2$ in $\Ad(K)$-invariant and irreducible subspaces, where $d_1 = 20$ and $d_2= 40$.  The only nonzero structural constant is $[122]=10$.  According to \eqref{scal}, the scalar curvature of a metric $g=(x,y)\in\mca^G$ is given by 
$$
\scalar(x,y) = \unm\left(\tfrac{20}{x}+\tfrac{40}{y}\right) -\tfrac{[122]}{4}\left(\tfrac{2}{x}+\tfrac{x}{y^2}\right) 
= \tfrac{5}{x}+\tfrac{20}{y}-\tfrac{5x}{2y^2}, 
$$
and if we set $x(y):=y^{-2}$, then $\scalar(y)=5y^2+\tfrac{20}{y}-\tfrac{5}{2y^4}$.  Since $\scalar'(1)=\scalar''(1)=0$ and $\scalar'''(1)=180\ne 0$, we obtain that $y=1$ is an inflection point of $\scalar(y)$ and hence $g_{\kil}$ is not a local maximum of $\scalar|_{\mca^G_1}$.

\section{Neutrally stable symmetric metrics}\label{NS}

In this section, we consider the only two cases among the neutrally stable irreducible symmetric spaces listed in \eqref{nsss} which admits a presentation $M=G/K$ such that $\dim{\mca_1^G}>0$.  We give in both cases an explicit curve $g(t)\in\mca_1^G$ through the symmetric metric $g_0\in\mca_1^G$ such that $\scalar(g(t))$ has an inflection point, in order to show that $g_0$ is not a local maximum of $\scalar|_{\mca_1^G}$ and so neither of $\widetilde{\scalar}|_{\cca_1}$.  This provides a proof for Theorem \ref{main}.

\subsection{Case $M=\SU(n)$, $n\geq 3$}
For $H:=\SU(n-1)\subset\SU(n)$, we consider the usual block matrix orthogonal decomposition in $\Ad(H)$-invariant and irreducible subspaces
$$
\sug(n)=\pg_1\oplus\pg_2\oplus\pg_3, \qquad \pg_1:=\sug(n-1), \quad \pg_2:=\CC^{n-1}, \quad\pg_3:=\RR A_0,
$$
where $A_0\in\sug(n)$ is the unique (up to scaling) diagonal matrix commuting with $H$, so $d_1= (n-1)^2-1$, $d_2 = 2(n-1)$ and $d_3 = 1$.  Using a well-known formula for the Einstein constant (see e.g.\ \cite[(8)]{stab-dos}), it is easy to see that the only nonzero structural constants are 
$$
[111] = (n-1)(n-2), \quad [122] = n-2, \quad [223]=1.  
$$
It follows from \eqref{scal} that for any metric $g=(x,y,z)$ such that $x^{d_1} y^{d_2}z=1$, the scalar curvature of $g$ is given by 
\begin{align*}
\scalar(x,y) =& \tfrac{1}{2}\left(\tfrac{d_1}{x}+\tfrac{d_2}{y}+\tfrac{d_3}{z}\right)-\tfrac{1}{4}\left([111]\tfrac{1}{x}+ [122]\left(\tfrac{x}{y^2} + \tfrac{2}{x}\right) + [223]\left(\tfrac{2}{z} + \tfrac{z}{y^2}\right)\right) \\ 
=& \tfrac{n^2-3n+2}{4x}+ \tfrac{n-1}{y} - \tfrac{(n-2)x}{4y^2} - \tfrac{1}{4}x^{-(n-1)^2+1} y^{-2n}.
\end{align*}
It can be seen that the kernel of the Hessian of this function at $(1,1)$ is generated by $(-\tfrac{2}{n-2},1)$, so we propose the curve of metrics $g(y)$ obtained by setting $x(y):=-\tfrac{2}{n-2} y + \tfrac{n}{n-2}$, which satisfies $g(1)=(1,1)=g_{\kil}$ and $g'(1)=(-\tfrac{2}{n-2},1)$.  A straightforward computation now gives 
\begin{align*}
\scalar(y) =& \tfrac{1}{4(n-2)y^2 x(y)} \left(y^2 n^2(n-5)+4n^2 y - n^2 + x(y)^{-n(n-2)}\left(2 y^{-2n+3} - n y^{-2n+2}\right) \right), \\ 
\scalar'(y) =& \tfrac{n^2}{2(n-2)^2 y^3 x(y)^2} \left(y^3 (n-5) + 8 y^2 - y (2n+3) + n
  + \tfrac{2 y^2 - (n+2) y + n}{x(y)^{n(n-2)}y^{2n-2}}\right), \\ 
\scalar''(y) =& \tfrac{n^2}{2(n-2)^3 y^4 x(y)^3} \left( 4(n-5) y^4 + 48 y^3 - 24(n+1) y^2  + 4n(4+n) y -3n^2 \right.\\ 
&+ x(y)^{-n(n-2)}y^{2-2n}(4(n^2+1) y^3 - 2(n^3+4n^2+n+4) y^2 \\ 
&\left.+ 2n (2n^2+2n+3) y -n^2 (2n+1)) \right).
\end{align*}
Thus 
$$
\scalar'(1)=0,  \qquad \scalar''(1)=0, \qquad \scalar'''(1)= \tfrac{n^2(n^3-5n^2+8n-4)}{(n-2)^4}=
\tfrac{n^2(n-1)}{(n-2)^2}\ne 0,
$$
that is, $y=1$ is an inflection point of $\scalar(y)$.  This implies that $g_{\kil}$ is not a local maximum of $\scalar|_{\mca^G_1}$.  In the case $n=3$, the curve $g(y)$ coincides with the one given in \cite[pp.1140]{Jns}.

\subsection{Case $M=\SU(2n)/\Spe(n)$, $n\geq 2$}
By considering the alternative presentation $M=G/K=\SU(2n-1)/\Spe(n-1)$, one obtains a reductive decomposition $\ggo=\kg\oplus\pg$ such that $\pg$ decomposes as $\pg=\pg_1\oplus\pg_2\oplus\pg_3$ in $\Ad(K)$-invariant and irreducible subspaces, with $d_1 = (2n-1)(n-2)$, $d_2= 4(n-1)$ and $d_3= 1$.

The scalar curvature of a metric $g=(x,y,z)\in\mca^G$ such that $x^{d_1} y^{d_2}z=1$ was computed in \cite{BhmWngZll}: 
\begin{align*}
\scalar(x,y) 
=& (2n-1)\left(\tfrac{4(n-1)(n-2)}{x} + \tfrac{8(n-1)}{y} - (n-2)\tfrac{x}{y^2} - x^{-(2n-1)(n-2)} y^{-4n+2}\right), 
\end{align*}
where it is also shown that the symmetric metric is given, up to scaling, by $g_0=(a,\tfrac{a}{2},\tfrac{n}{2n-1}a)$ for any $a>0$.  In order to get the right volume, we must set 
$$
a:= \left(\tfrac{16n}{(2n-1) 16^{n}}\right)^{\tfrac{1}{n+1-2n^2}}.
$$  
It can be seen that the kernel of the Hessian of this function at $(a,\tfrac{a}{2})$ is generated by $(1,-\tfrac{n-2}{4})$, so we take the curve of metrics $g(x)\in\mca^G_1$ obtained by setting $y(x):=-\tfrac{n-2}{4} x + \tfrac{n}{4}a$, which satisfies $g(1)=(1,\unm)=g_0$ and $g'(1)=(1,-\tfrac{n-2}{4})$.  One therefore obtains that
\begin{align*}
\scalar(x) =& (2n-1)\left(\tfrac{4(n-1)(n-2)}{x} + \tfrac{8(n-1)}{y(x)} - (n-2)\tfrac{x}{y(x)^2}- x^{-(2n-1)(n-2)} y(x)^{-4n+2}\right) \\ 
=& (2n-1)(p(x) - q(x)),
\end{align*}
where
$$
p(x):= \tfrac{4(n-1)(n-2)}{x} + \tfrac{8(n-1)}{y(x)} - (n-2)\tfrac{x}{y(x)^2}, \qquad
q(x):=  \tfrac{1}{x^{d_1} y(x)^{d_2+2}}.  
$$
Since
\begin{align*}
p'(x)=&-\tfrac{4(n-1)(n-2)}{x^2} - \tfrac{8(n-1)y'(x)}{y^2(x)} - \tfrac{(n-2)(y(x)^2-2xy(x)y'(x))}{y(x)^4} \\ 
   =& -\tfrac{4(n-1)(n-2)}{x^2} + \tfrac{2(n-1)(n-2)}{y^2(x)} - \tfrac{(n-2)(2y(x) +(n-2)x)}{2 y(x)^3}, \\ 
p''(x) =& \tfrac{8(n-1)(n-2)}{x^3} - \tfrac{(n-1)(n-2)^2}{y^3(x)} + \tfrac{(n-2)^2(3x-y(x))}{2 y(x)^4}, \\
q'(x)=& -d_1 x^{-d_1-1} y(x)^{-d_2-2} - (d_2+2)x^{-d_1} y(x)^{-d_2-3}y'(x) \\ 
 =& - x^{-(2n-1)(n-2)} y(x)^{-(4n-2)} \left(\tfrac{(2n-1)(n-2)}{x} - \tfrac{(n-2)(2n-1)}{2 y(x)}\right)  \\ 
q''(x) =& d_1(d_1+1) x^{-d_1-2} y(x)^{-d_2-2} 
+ 2 d_1 (d_2+2) x^{-d_1-1} y(x)^{-d_2-3}y'(x) \\ 
&+ (d_2+2)(d_2+3)x^{-d_1} y(x)^{-d_2-4}y'(x)^2,
\end{align*}
a straightforward computation gives that
$$ 
p'(a)=0,\quad p''(a)=\tfrac{2(n-2)n^2}{a^3}, \quad p'''(a)=-\tfrac{6(n-2)n^2}{a^4},
$$
and 
$$
q'(a)=0,\quad q''(a)=\tfrac{n (2n-1)(n-2) 16^n}{8 a^{2n^2-n+2}}, \quad q'''(a)=\tfrac{2n^2(n-2)(2n^2-9n+4)}{(2n-1)a^4}.
$$
Thus $\scalar'(a)=0$, $\scalar''(a)=0$ and $\scalar'''(a)\ne 0$, that is, $x=a$ is an inflection point of $\scalar(x)$, from which we deduce that $g_0$ is not a local maximum of $\scalar|_{\mca^G_1}$.    

This was proved in \cite[Example 6.7]{BhmWngZll} by using different methods.

\section{$G$-neutrally stable full flag manifold}\label{GNS}

We study in this section the only known $G$-neutrally stable Einstein metric with $\dim{\mca_1^G}\geq 2$.  It was found in \cite[Section 5]{stab}.  

Consider the standard metric $g_{\kil}$ on the homogeneous spaces $\SO(2n)/T^n$, $n\geq 4$, where $T^n$ is the usual maximal torus of $\SO(2n)$, which is known to be Einstein (see \cite{WngZll2}).  The standard block matrix reductive decomposition is given by 
$$
\ggo=\kg\oplus\pg, \qquad\pg=\pg_{12}\oplus\pg_{13}\oplus\dots\oplus\pg_{(n-1)n},
$$
where every block $\pg_{ij}=\pg_{ji}$ (note that always $i\ne j$) has dimension $4$ and they are all $\Ad(T^n)$-invariant but not irreducible (see \cite[Section 5]{stab-tres}).  All the above triples produces the same nonzero structural constant $\tfrac{2}{n-1}$.

We consider the metric $g=(x,\dots,x,y,\dots,y)$, where the value $x$ corresponds to $\pg_{1j}$, $j=2,\dots,n$ and $y$ otherwise.  It follows from \eqref{scal} that 
$$  
\scalar(x,y)= \tfrac{2(n-1)}{x} + \tfrac{(n-1)(n-2)}{2y} - \tfrac{n-2}{2}\tfrac{y}{x^2},
$$
and if we set $x^{n-1}y^{\tfrac{(n-1)(n-2)}{2}}=1$, then it is straightforward to see that
\begin{align*}
\scalar(x)=& \tfrac{4(n-1) + (n^2-3n+2)x^{\tfrac{n}{n-2}} + (2-n) x^{-\tfrac{n}{n-2}}}{2x}, \\ 
\scalar'(x)=& -\tfrac{2(n-1)}{x^2}+\tfrac{n-1}{x} x^{\tfrac{2}{n-2}} + (n-1)x^{-\tfrac{2}{n-2}-3},\\ 
\scalar''(x)=& \tfrac{4(n-1)}{x^3}+\tfrac{(n-1)(4-n)}{n-2} x^{\tfrac{2}{n-2}-2} - \left(\tfrac{2}{n-2} + 5 +3(n-2)\right) x^{-\tfrac{2}{n-2}-4}.  
\end{align*}
Thus $\scalar'(1)=\scalar''(1)=0$ and $\scalar'''(1)=2 n^2 \tfrac{(n-1)}{(n-2)^2}\ne 0$, so $x=1$ is an inflection point of $\scalar(x)$ and $g_{\kil}$ is not a local maximum of $\scalar|_{\mca^G_1}$.


\begin{thebibliography}{MMM}

%\bibitem[A]{Alk} {\sc D. Alekseevskii}, Conjugacy of polar factorizations of Lie groups, {\it Mat. Sb.} {\bf 84} (1971), 14-26; {\it English translation}: {\it Math. USSR-Sb.} {\bf 13} (1971), 12-24.

%\bibitem[APZ]{ArrPlmZll} {\sc R. Arroyo, A. Pulemotov, W. Ziller}, The prescribed Ricci curvature problem for naturally reductive metrics on compact Lie groups, preprint 2020 (arXiv).
%
%\bibitem[AGP]{ArrGldPlm} {\sc R. Arroyo, M. Gould, A. Pulemotov}, The prescribed Ricci curvature problem for naturally reductive metrics on non-compact simple Lie groups, preprint 2020 (arXiv).

\bibitem[BHMW]{BttHllMrpWld} {\sc W. Batat, S. Hall, T. Murphy, J. Waldron}, Rigidity of $\SU_n$-type symmetric spaces, preprint 2021 (arXiv).  

\bibitem[B]{Bss} {\sc A. Besse}, Einstein manifolds, {\it Ergeb. Math.} {\bf 10} (1987), Springer-Verlag, Berlin-Heidelberg.

%\bibitem[BCR]{BchCstRoy} {\sc J. Bochnak, M. Coste, M.-F. Roy}, Real Algebraic Geometry, {\it Ergeb. Math.} {\bf 36} (1998), Springer-Verlag.
%
%\bibitem[Bo]{Bhm} {\sc C. B\"ohm}, Unstable Einstein metrics, {\it Math. Zeit.} {\bf 250} (2005), 279-286.  
%
%\bibitem[BK]{BhmKrr} {\sc C. B\"ohm, M. Kerr}, Low-dimensional homogeneous Einstein manifolds, {\it Trans. Amer. Math. Soc.} {\bf 358} (2005), 1455-1468. 
%
%\bibitem[BL]{BhmLfn} {\sc C. B\"ohm, R. Lafuente}, Real geometric invariant theory, {\it Diff. Geom. in the Large} (2020), Cambridge Univ. Press, in press.

\bibitem[BWZ]{BhmWngZll} {\sc C. B\"ohm, M.Y. Wang, W. Ziller}, A variational approach for compact homogeneous Einstein manifolds, {\it Geom. Funct. Anal.} {\bf 14} (2004), 681-733.

%\bibitem[Br]{Brd} {\sc G. Bredon}, Introduction to compact transformation groups, {\it Pure and Applied
%Mathematics} {\bf 46} (1972), Academic Press.  

%\bibitem[Bu]{Btt} {\sc T. Buttsworth}, The prescribed Ricci curvature problem on three-dimensional unimodular Lie groups, {\it Math. Nachr.} {\bf 292} (2019), 747-759.

%\bibitem[BH]{BttHll} {\sc T. Buttsworth, M. Hallgren}, Local stability of Einstein metrics under the Ricci iteration, preprint 2019 (arXiv).

%\bibitem[BP]{BttPlm} {\sc T. Buttsworth, A. Pulemotov}, The prescribed Ricci curvature problem for homogeneous metrics, {\it Differential geometry in the large}, Cambridge University Press, in press.
%
%\bibitem[BPRZ]{BttPlmRbnZll} {\sc T. Buttsworth, A. Pulemotov, Y.A. Rubinstein, W. Ziller}, On the Ricci iteration for homogeneous metrics on spheres and projective spaces, preprint 2018 (arXiv).

%\bibitem[CHI]{CaoHmlIlm} {\sc Huai-Dong Cao, R.S. Hamilton, T. Ilmanen}, Gaussian densities and
%stability for some Ricci solitons, preprint 2004 (arXiv).

\bibitem[CH]{CaoHe} {\sc Huai-Dong Cao, Chenxu He}, Linear stability of Perelman's $\nu$-entropy on symmetric spaces of compact type, {\it J. reine angew. Math.}, {\bf 709} (2015), 229-246.

%\bibitem[DZ]{DtrZll} {\sc J. D'Atri, W. Ziller}, Naturally reductive metrics and Einstein metrics on compact lie groups, {\it Mem. Amer. Math. Soc.} {\bf 215} (1979).

%\bibitem[De]{Dly} {\sc E. Delay}, ....
%
%\bibitem[D]{Dtr} {\sc D. DeTurck}, Prescribing positive Ricci curvature on compact manifolds, {\it Rend. Sem. Mat. Univ. Politec. Torino} {\bf 43} (1985) 357-369.
%
%\bibitem[DK]{DtrKso} {\sc D. DeTurck, N. Koiso}, Uniqueness and nonexistence of metrics with prescribed Ricci curvature, {\it Ann. Inst. Poincar\'e Anal. Non Lin\'eaire 1} (1984) 351-359.
%

\bibitem[DK]{DckKrr} {\sc W. Dickinson, M. Kerr}, The geometry of compact homogeneous spaces with two isotropy summands, {\it Ann. Global Anal. Geom.} {\bf 34} (2008), 329-350. 

%\bibitem[DLM]{DttLtMtl} {\sc I. Dotti, M. L. Leite, R. Miatello}, Negative Ricci curvature on complex semisimple Lie grous, {\it Geom. Dedicata} {\bf 17} (1984), 207-218.
%

\bibitem[GG]{GsqGld} {\sc J. Gasqui, H. Goldschmidt}, Radon transforms and spectral rigidity on the
complex quadrics and the real Grassmannians of rank two, {\it J. Reine Angew. Math.} {\bf 480} (1996), 1-69.

%\bibitem[GR]{GdsRyl} {\sc C. Godsil, G. Royle}, Algebraic graph theory, {\it GTM} {\bf 207} (2001), Springer.   
%
%\bibitem[GS]{GttSmm} {\sc S. Goette, U. Semmelmann}, Scalar curvature estimates for compact symmetric spaces, {\it Diff. Geom. App} {\bf 16} (2002), 65-78.
%
%\bibitem[G]{Grd} {\sc C. Gordon}, Naturally reductive homogeneous Riemannian manifolds, {\it Can. J. Math.} {\bf 37} (1985) 467-487.

%\bibitem[G1]{Grm2} {\sc M. Gromov}, Positive curvature, macroscopic dimension, spectral gaps and higher signatures, in: S. Gindikin, J. Lepowski and R. L. Wilson.eds., Functional Analysis on the Eve of the 21st Century, Vol. II, Progress in Mathematics Vol. 132 (1996), 1?213.
%
%\bibitem[G2]{Grm} {\sc M. Gromov}, Four Lectures on Scalar Curvature, preprint 2019 (arXiv).

%\bibitem[H]{Hbr} {\sc J. Heber}, Noncompact homogeneous Einstein spaces, {\it Invent. math}. {\bf 133} (1998), 279-352.

%\bibitem[Ha]{Hml} {\sc R. Hamilton}, The Ricci curvature equation, in: Seminar on nonlinear partial differential equations (S.-S. Chern, ed.), Math. Sci. Res. Inst. Publ. 2, Springer-Verlag, New York, 1984, 47-72.

\bibitem[H]{He} {\sc Chenxu He}, Cohomogeneity one manifolds with a small family of invariant metrics, {\it Geom. Dedicata} {\bf 157}, 41-90.

\bibitem[J]{Jns} {\sc G. Jensen}, The Scalar Curvature of Left-Invariant Riemannian Metrics, {\it Indiana Math. J.} {\bf 20} (1971), 1125-1144.

\bibitem[Ke]{Krr} {\sc M. Kerr}, New Examples of Homogeneous Einstein Metrics, {\it Michigan Math. J.} {\bf 45} (1998), 115-134.

\bibitem[K]{Kso} {\sc N. Koiso}, Rigidity and stability of Einstein metrics: the case of compact symmetric spaces, {\it Osaka J. Math.} {\bf 17} (1980), 51-73.

%\bibitem[Ko]{Kst2} {\sc B. Kostant}, On differential geometry and homogeneous spaces I and II, {\it Proc. Nat. Acad. Sci. U.S.A.} {\bf 42} (1956), 258-261and 354-357.

\bibitem[Kr1]{Krn2} {\sc K. Kr\"oncke}, Stability of Einstein Manifolds, {\it Ph.D. thesis} (2013), Universit\"at Potsdam.  

\bibitem[Kr2]{Krn} {\sc K. Kr\"oncke}, Stability and instability of Ricci solitons, {\it Calc. Var. PDE.}  {\bf 53} (2015), 265-287.

%\bibitem[LL]{alek} {\sc R. Lafuente, J. Lauret}, Structure of homogeneous Ricci solitons and the Alekseevskii conjecture, {\it J. Diff. Geom.} {\bf 98} (2014), 315-347.

\bibitem[LL]{stab} {\sc E.A. Lauret, J. Lauret}, The stability of standard homogeneous Einstein manifolds, preprint 2021 (arXiv).  

%\bibitem[L1]{natred}  {\sc J. Lauret}, Naturally reductive homogeneous structures on $2$-step nilpotent Lie groups, {\it Rev. Un. Mat. Argentina} {\bf 41} (1998), 15-23.

%\bibitem[L2]{manus} {\sc J. Lauret}, Homogeneous nilmanifolds attached to representations of compact Lie groups, {\it Manusc. Math.} {\bf 99} (1999), 287-309.
%
%\bibitem[L3]{homRF}  {\sc J. Lauret}, Ricci flow of homogeneous manifolds, {\it Math. Z.} {\bf 274} (2013), 373-403.
%
%\bibitem[L4]{sol-HS}  {\sc J. Lauret}, The search for solitons on homogeneous spaces, {\it Abel Symposia} (2019), Springer, in press.

%\bibitem[LW1]{RNder} {\sc J. Lauret, C.E. Will}, On Ricci negative Lie groups, {\it Abel Symposia} (2019), Springer, in press.

\bibitem[L]{stab-tres} {\sc J. Lauret}, On the stability of homogeneous Einstein manifolds, {\it Asian J. Math.}, in press (arXiv).

%\bibitem[LW1]{PRP} {\sc J. Lauret, C.E. Will}, Prescribing Ricci curvature on homogeneous manifolds, {\it J. reine angew. Math.}, in press (arXiv).  

\bibitem[LW]{stab-dos} {\sc J. Lauret, C.E. Will}, On the stability of homogeneous Einstein manifolds II, preprint 2021 (arXiv).

%\bibitem[N]{Nkn} {\sc Y.G. Nikonorov}, The scalar curvature functional and homogeneous Einsteinian metrics on Lie groups, {\it Siberian Math. J.} {\bf 39} (1998) 504-509.
%

\bibitem[N1]{Nkn} {\sc Yu.G. Nikonorov}, On a characterization of critical points of the scalar curvature functional (Russian), {\it Tr. Rubtsovsk. Ind. Inst.} {\bf 7} (2000), 211-217.  {\it Translation}: arXiv:2112.00993.  This article consists of \cite[Section 1.2, pp. 25-34]{Nkn2}.  

\bibitem[N2]{Nkn2} {\sc Yu.G. Nikonorov}, Analytical methods in the theory of homogeneous Einstein manifolds (Russian), Dissertation (thesis) for the degree of Doctor of Physical and Mathematical Sciences, Rubtsovsk, 2002.

%\bibitem[NRS]{NknRdnSlv} {\sc Y.G. Nikonorov, E.D. Rodionov, V.V. Slavskii}, Geometry of homogeneous Riemannian manifolds, {\it J. Math. Sci. (N.Y.)} {\bf 146} (2007) 6313-6390.
%
%\bibitem[Pe]{Prl} {\sc G. Perelman}, The entropy formula for the Ricci flow and its geometric applications, preprint 2002 (arXiv).

%\bibitem[PR]{PlmRbn} {\sc A. Pulemotov, Y.A. Rubinstein}, Ricci iteration on homogeneous spaces, {\it Trans. AMS} {\bf 371} (2019), 6257-6287.

\bibitem[SW]{SmmWng} {\sc U. Semmelmann, G. Weingart}, Stability of Compact Symmetric Spaces,  {\it J. Geom. Anal.} {\bf 32} (2022), 137.

\bibitem[S]{Sch} {\sc P. Schwahn}, Stability of Einstein metrics on symmetric spaces of compact type, {\it Ann. Global Anal. Geom.} (2021), in press.

\bibitem[SSW]{SchSmmWng} {\sc P. Schwahn, U. Semmelmann, G. Weingart}, Stability of the Non-Symmetric Space $E_7/PSO(8)$, preprint 2022 (arXiv).  

%\bibitem[S]{Str} {\sc R. Storm}, The Classification of $7$- and $8$-dimensional Naturally Reductive Spaces, {\it Canadian. J. Math.}, in press.  
%
%\bibitem[SD]{SunDai} {\sc Yukai Sun, Xianzhe Dai}, Extremality of Bi-invariant Metrics on Lie Groups and
%Homogeneous Spaces, preprint 2020 (arXiv).

%\bibitem[W]{Wng} {\sc M. Y. Wang}, Einstein Metrics from Symmetry and Bundle Constructions: A Sequel, {\it Adv. Lect. Math.} {\bf 22} (2012), 253-309, International Press. 

\bibitem[WW]{WngWng} {\sc Changliang Wang, M. Y. Wang}, Instability of some Riemannian manifolds with real Killing spinors, {\it Comm. Anal. Geom.} (2021), in press.  

\bibitem[WZ1]{WngZll2} {\sc M. Y. Wang, W. Ziller}, On normal homogeneous Einstein manifolds, {\it Ann.
Sci. \'Ecole Norm. Sup.} {\bf 18} (1985), 563-633.

\bibitem[WZ2]{WngZll} {\sc M.Y. Wang, W. Ziller}, Existence and nonexistence of homogeneous Einstein metrics, {\it Invent. Math.} {\bf 84} (1986), 177-194.
\end{thebibliography}
\end{document}